\documentclass{article}
\usepackage{graphicx}

\title{A Hybrid Framework for Healing Semigroups with Machine Learning}
\author{
  Sarayu Sirikonda\thanks{Mountain House High School}
  \and
  Jasper van de Kreeke\thanks{University of California, Berkeley}
}
\date{}

\usepackage[table,dvipsnames]{xcolor}
\usepackage{tikz, float, tocloft, pgfplots, booktabs, tabularx, siunitx, sectsty, comment, enumitem, cite}
\usepackage{amsmath,amssymb}
\usetikzlibrary{arrows.meta, positioning, shapes.multipart, matrix, fit}
\pgfplotsset{compat=1.17}

\usepackage{hyperref}
\hypersetup{colorlinks=true,linkcolor=blue,citecolor=blue,urlcolor=blue}

\begin{document}

\setlength{\cftsecindent}{0.1pt} 
\setlength{\cftsubsecindent}{1em}
\setlength{\cftsecnumwidth}{2em}
\setlength{\cftsubsecnumwidth}{2.5em}
\setlength{\cftbeforesecskip}{0pt}
\sectionfont{\centering}

\maketitle
\begin{abstract}
In this paper, we propose a hybrid framework that heals corrupted finite semigroups, combining deterministic repair strategies with Machine Learning using a Random Forest Classifier. Corruption in these tables breaks associativity and invalidates the algebraic structure. Deterministic methods work for small cardinality $n$ and low corruption but degrade rapidly. Our experiments, carried out on Mace4-generated data sets, demonstrate that our hybrid framework achieves higher healing rates than deterministic-only and ML-only baselines. At a corruption percentage of $p=15\%$, our framework healed 95\% of semigroups up to cardinality $n=6$ and 60\% at $n=10$.

\end{abstract}

\setcounter{tocdepth}{1}
\tableofcontents 

\section{Introduction}
\paragraph{Background.}
Semigroups are among the simplest algebraic structures. However, they appear across mathematics and applications, serving as a foundational structure in computer science, biology, and physics \cite{li2022semigroup,cuturi2005semigroup,feng2017sgdpca,assiry2023roughness}. 

In this work, we focus on \textit{finite} semigroups. A finite semigroup of size $n$ can be represented as an $n \times n$ table (also known as a Cayley table), where the entry in row $i$ and column $j$ holds the result of the operation $i \cdot j$. Each element in the table is an integer in the set $\{0, 1, \dots, n-1\}$.

\paragraph{Motivation.}
While much of algebra focuses on studying perfectly defined structures, real-world computations often yield \emph{broken} or corrupted structures. Studying how to repair broken algebraic structures is an important step toward developing robust methods and algorithms, and semigroups provide a natural starting point for this investigation because of their limited structure.

\paragraph{Deterministic Healing.}
As a baseline, we designed two deterministic repair strategies. While simple, they achieved at most 47\% associativity at cardinality $n=3$, with performance deteriorating further as $n$ increased (see Section 5). 

Although extensions such as backtracking improved associativity, they failed to preserve per-cell accuracy, producing repaired structures that diverged from the original.

These limitations demonstrate that purely deterministic methods cannot reliably scale beyond very small semigroups and motivates the development of our hybrid framework. 

\section{Related Work}

Even for small sizes, the number of possible $n \times n$ multiplication tables grows as $n^{n^2}$, so brute-force search is infeasible. This motivates algorithmic and heuristic methods that construct or analyze semigroups without enumerating all tables. The Froidure–Pin algorithm remains a foundational method to compute the structure of a finite semigroup from generators, producing Cayley graphs and rewriting data used by many later systems. Counting results underscore the combinatorial explosion: Distler, Jefferson, Kelsey, and Kotthoff proved that there are $12{,}418{,}001{,}077{,}381{,}302{,}684$ non-equivalent semigroups of order 10 \cite{distler2012semigroups10}; related enumerations for monoids of orders 8–10 appeared earlier. These results explain why exhaustive methods quickly become infeasible and motivate learning-based surrogates or priors \cite{williamson2023deeplearning,koehler2024multimodal}.

More recently, researchers have started exploring the use of machine learning for algebraic structures. The most directly relevant effort is due to Balzin and Shminke (2021), who introduced an autoencoder framework to reconstruct semigroup tables with missing entries. Their method demonstrated that partial information is sufficient to recover associative tables in many cases \cite{balzin2021neural}.

A separate line of research encodes the semigroup property in learning systems for \textit{operators} rather than finite algebras: learning generators of Markov semigroups via resolvent/Laplace methods \cite{kostic2024laplace}, analyzing diffusion-based GNNs through operator-semigroup ergodicity \cite{zhao2024oversmoothing}, embedding semigroup composition laws into neural operators for ODE/PDE flow maps \cite{chen2025due}, and proving contraction/exponential convergence of iterative transport algorithms viewed as semigroups \cite{akyildiz2025sinkhorn}. These results show that respecting semigroup structure can improve stability, sample efficiency, and provable guarantees—but they concern continuous/operator settings, not discrete Cayley tables \cite{hillebrecht2025pinnserror}.

\section{Preliminaries}

\paragraph{Semigroups.} A \emph{semigroup} is a pair $(S,\cdot)$ where $S$ is a finite set of size $n$ and $\cdot:S\times S \to S$ is a binary operation that is \emph{associative}:
\[
(a \cdot b)\cdot c \;=\; a \cdot (b \cdot c) \quad \forall a,b,c \in S.
\]

\paragraph{Example.} 
The integers modulo $3$ with addition form a semigroup $(\{0,1,2\},+)$, with Cayley table.
\[
\begin{array}{c|ccc}
+ & 0 & 1 & 2 \\
\hline
0 & 0 & 1 & 2 \\
1 & 1 & 2 & 0 \\
2 & 2 & 0 & 1
\end{array}
\]

\paragraph{Isomorphism.} 
Two semigroups $(S,\cdot)$ and $(S',\ast)$ of the same size are \emph{isomorphic} if there exists a bijection $\varphi:S\to S'$ such that $\varphi(a\cdot b)=\varphi(a)\ast\varphi(b)$ for all $a,b\in S$. Equivalently, if we label $S=\{0,\dots,n-1\}$, an isomorphism is a permutation $\pi\in S_n$ such that for the Cayley tables $T,T'$:
\[
T'\big[\pi(i),\pi(j)\big] = \pi(T[i,j]) \quad \forall i,j.
\]
This means many different labeled Cayley tables represent the same algebraic structure; when counting semigroups we report distinct classes up to isomorphism. This distinction is important, since the number of distinct semigroups grows extremely quickly as $n$ increases.

\paragraph{Growth.} 
Classical enumerations illustrate this combinatorial explosion:
\[
\begin{array}{c|cccccccc}
n & 1 & 2 & 3 & 4 & 5 & 6 & 7 & 8 \\ \hline
\#\text{ semigroups} & 1 & 4 & 18 & 126 & 1,160 & 15,973 & 836,021 & 1,843,120,128
\end{array}
\]
By $n=9$, the number already exceeds billions, underscoring both the richness of the space and the challenge for any data-driven approach.

\subsection*{Definition of Healing}
Given a corrupted Cayley table $T$, the goal of \emph{healing} is to produce a repaired table $\hat{T}$ that restores the semigroup structure. For this paper, a healed table must satisfy two requirements:
\begin{enumerate}
    \item \textbf{Global associativity:} $(a \cdot b) \cdot c = a \cdot (b \cdot c)$ for all $a,b,c \in S$.
    \item \textbf{Local fidelity:} $\hat{T}$ should remain close to the original uncorrupted table $T^\ast$, measured by per-cell accuracy.
\end{enumerate}

\section{Experimental Setup}

\begin{figure}
    \centering
    \includegraphics[width=0.9\linewidth]{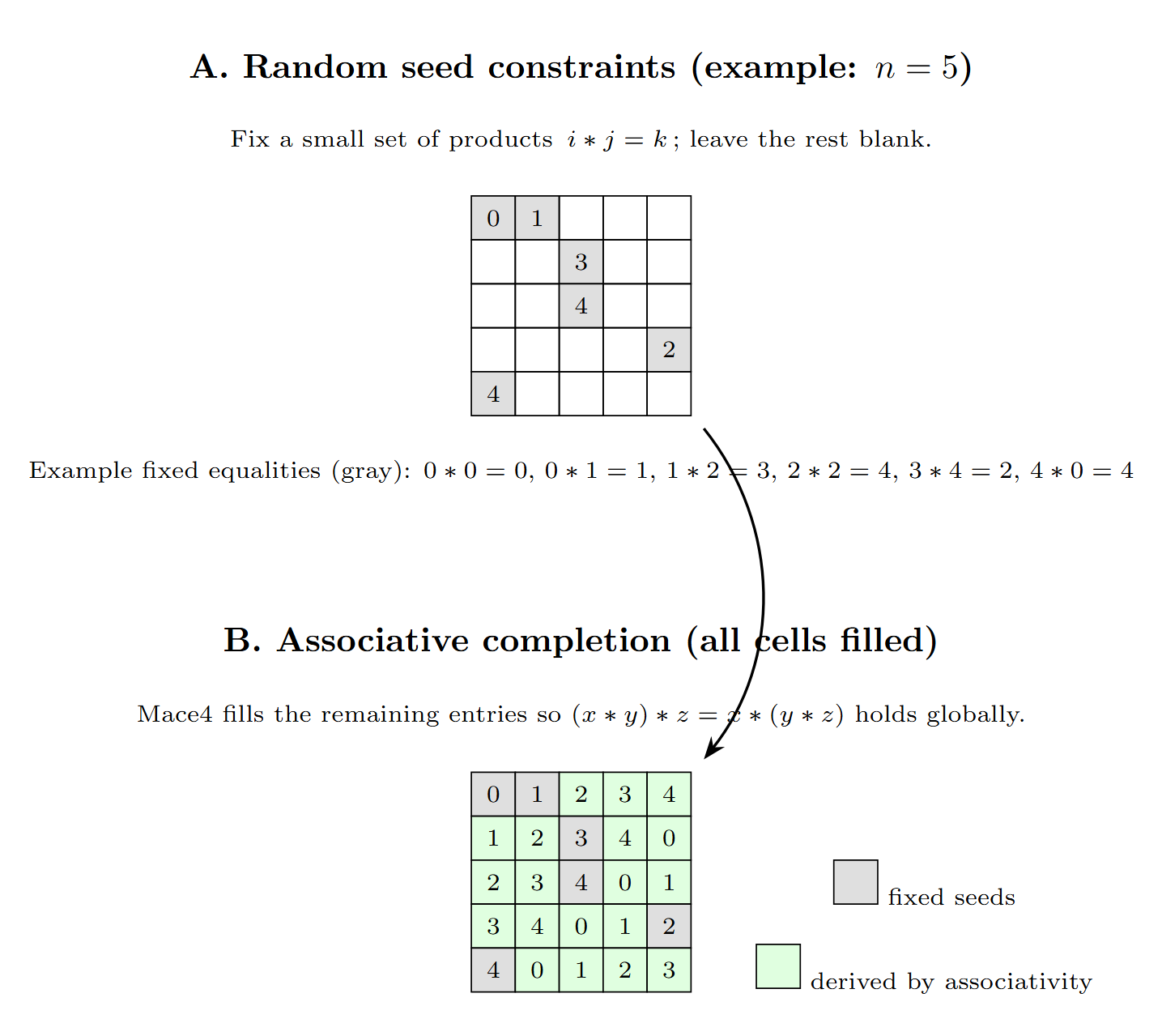}
    \caption{\textbf{Mace4 generation model (5$\times$5).} (A) A handful of seed entries (gray) are fixed. (B) The table is deterministically completed (green) so associativity holds for all triples.}
    \label{fig:mace4}
\end{figure}

\subsection{Dataset Generation}
\paragraph{Brute Force.} 
One of the methods we used was a "brute force" approach, in which table entries were filled randomly with values $0$ to $n-1$, until full-table associativity was satisfied. This method proved to be highly inefficient: taking 4 minutes to produce just 10 associative tables of cardinality $n=4$. 

\paragraph{Mace4 Generation.} To address this limitation, we use the \texttt{Mace 4} program, a finite model builder \cite{mace4}. This tool generates hundreds of associative tables up to $n \approx 30$ within seconds, performing much better than brute force.

Once correct semigroup tables are generated, we select $p\%$ (corruption percentage) of entries uniformly at random and flip each to an incorrect value from $\{0,\dots,n-1\}\setminus\{T[i,j]\}$. This preserves closure but breaks associativity. The resulting pairs of $(\text{clean},\text{corrupt})$ tables form the dataset.

\renewcommand{\arraystretch}{1.2}
\[
\underbrace{\vcenter{\hbox{$
\begin{array}{c|ccc}
\cdot & 0 & 1 & 2\\ \hline
0 & 0 & 1 & 2\\
1 & 1 & 2 & 0\\
2 & 2 & 0 & 1\\
\end{array}
$}}}_{\text{Clean}}
\quad\Longrightarrow\quad
\underbrace{\vcenter{\hbox{$
\begin{array}{c|ccc}
\cdot & 0 & 1 & 2\\ \hline
0 & 0 & 1 & 2\\
1 & 1 & \color{red}{0} & 0\\
2 & 2 & 0 & \color{red}{2}\\
\end{array}
$}}}_{\text{Corrupted ($p{=}15\%$)}}
\]


\subsection{Trust Maps}
Our hypothesis was that corrupted values were more likely to participate in failed associativity checks. To quantify and test this, we define a \emph{trust map} over the table.

For each entry $T(i,j)$ we compute its trust score as
\[
\text{trust}(i,j) \;=\; 
\frac{\#\{\,k \;\mid\; T(T(i,j),k) \;=\; T(i,T(j,k)) \,\}}
     {n}.
\]
That is, we count the number of associativity checks involving $T[i,j]$ that
are satisfied, divided by the total number of such checks. Equivalently, the
trust score represents the empirical likelihood that $T[i,j]$ is consistent
with associativity. Each entry therefore receives a value between 0 and 1.

The resulting trust map is a matrix of the same dimensions as the Cayley
table. High-trust cells are likely to be uncorrupted, while low-trust cells
are likely corrupted.

\begin{figure}[H]
    \centering
    \includegraphics[width=0.8\linewidth]{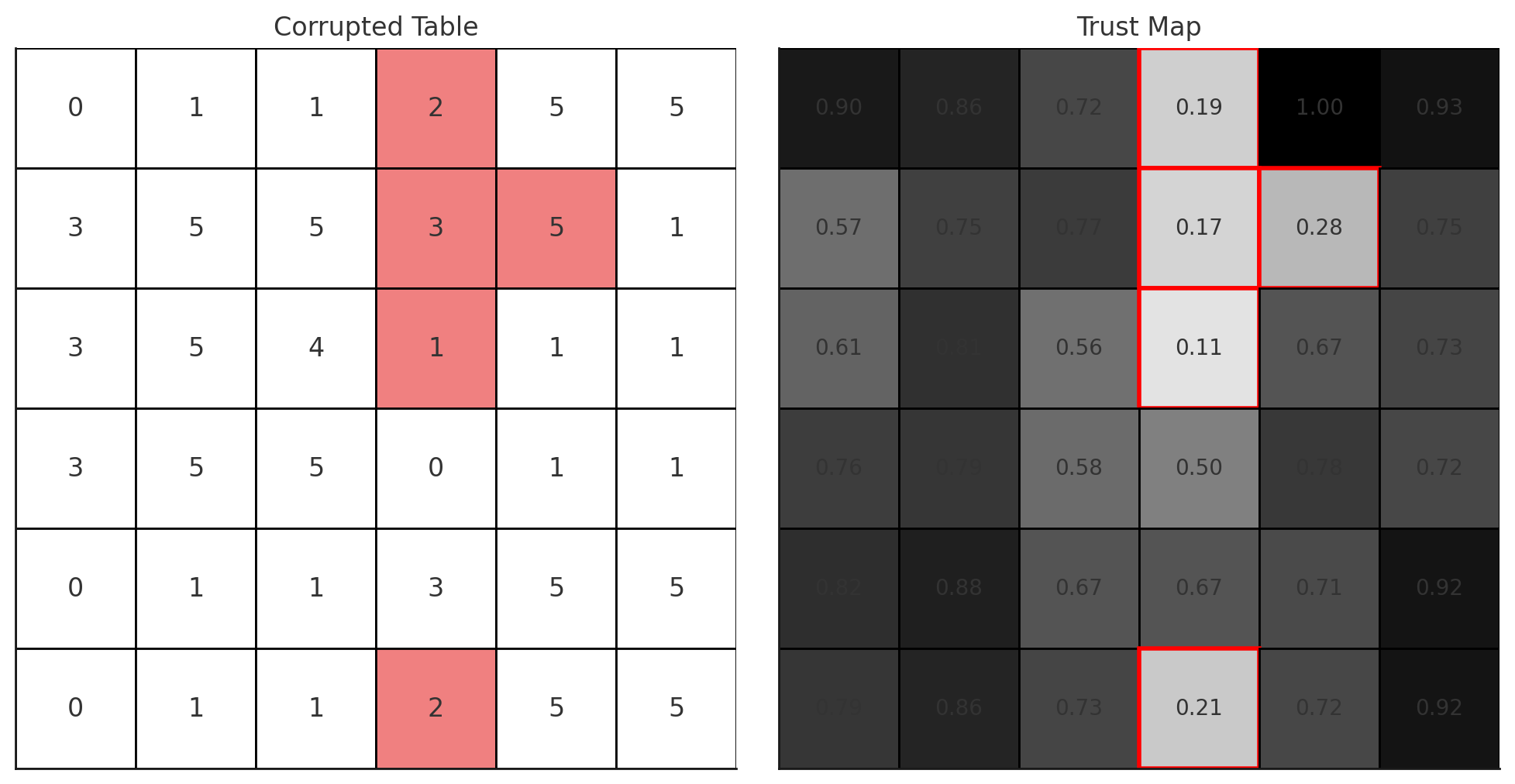}
    \caption{Trust map vs Corrupted table}
    \label{fig:trust-map}
\end{figure}

\textbf{Trust Map Accuracy.} To demonstrate the effectiveness of our trust map in identifying corrupted cells, we compare a corrupted semigroup table with its corresponding trust map (Figure \ref{fig:trust-map}). In the first table, the corrupted cells are highlighted in red. In the second table, we show the trust map values for each corresponding cell. Notably, the highlighted positions in the trust map align closely with the truly corrupted cells, demonstrating that the trust scoring mechanism effectively assigns lower scores to incorrect values. 

\subsection{Subsemigroups}
Directly healing a full semigroup table of higher cardinality proved to be almost impossible. To overcome this, we propose a layered healing strategy. We decompose the semigroup into smaller, overlapping \textit{subsemigroups}. Each is healed individually and then reconstructed back into the global structure, while ensuring that the final table remains associative. 

\paragraph{Local closure construction.}
For each triple $(i,j,k) \in  n^3$, we define the closure set: \[
    G(i,j,k) = \{i,j,k,i \cdot j,j \cdot k,(i \cdot j) \cdot k,i \cdot (j \cdot k)\}
\]

This is the smallest set of elements required to check associativity for $(i,j,k)$. We restrict attention to $|G| \in [2,5]$, since larger closures become redundant and smaller ones (size $1$) are trivial.

\paragraph{Validating Candidates.}
Not every $G(i,j,k)$ forms a valid subsemigroup. To guarantee closure we test \textit{every} product $a\cdot b$ and $a\cdot b \cdot c$ for $a,b,c \in G$:
\begin{itemize}
    \item If the result lies outside $G$, we reject this candidate.
    \item Otherwise, we keep it as a valid subsemigroup for repair.
\end{itemize}

This ensures each retained subsemigroup is algebraically meaningful and able to be repaired. 

\paragraph{Local Healing.}
Each surviving subsemigroup is reindexed (e.g., ${2,7,8}\mapsto{0,1,2}$) to simplify computation. On this smaller table:
\begin{enumerate}
    \item We construct a local trust map, measuring how reliable each entry (Section 3.2).
    \item For every triple, if $(i\cdot j)\cdot k \neq i\cdot (j\cdot k)$, we replace the side with lower trust by the other.
    \item The reindexed table is then mapped back to its global element labels.
\end{enumerate}

This enforces associativity within each local table.

\paragraph{Global merging.}
Since different subsemigroups overlap, multiple candidates may propose values for the same global entry. To resolve these conflicts, we assign a weight to each candidate: 
\[
    w = \frac{1}{|G|^2} \times \text{(ML Probability)} \times \text{(Trust Score)}
\]
and pick the candidate with the largest weight. The ML probability scores how likely a value is correct given learned patterns.

\paragraph{Justification.}
Every associativity check in the global table depends only on the small closure set $G(i,j,k)$. By repairing all such subtables, we guarantee local consistency for every triple. When these repaired subtables are reconstructed with weights, the result is a globally associative table.

\section{Deterministic Algorithm}

The deterministic repair procedure for a single entry $i\cdot j$ works by exploiting associativity. 
We enumerate all decompositions $i=i_1\cdot i_2$. 
For each, associativity gives
\[
(i_1\cdot i_2)\cdot j \;=\; i_1\cdot (i_2\cdot j).
\]
Each such identity yields a candidate value for $i\cdot j$. 
We tally the frequency of each candidate across all decompositions and set $i\cdot j$ to the most frequently supported value. 
This “majority vote” ensures the chosen entry is maximally consistent with associativity.

\subsection{Deterministic Results}
\begin{figure}[h]
    \centering
    \begin{minipage}[b]{0.48\textwidth}
        \centering
        \includegraphics[width=\linewidth]{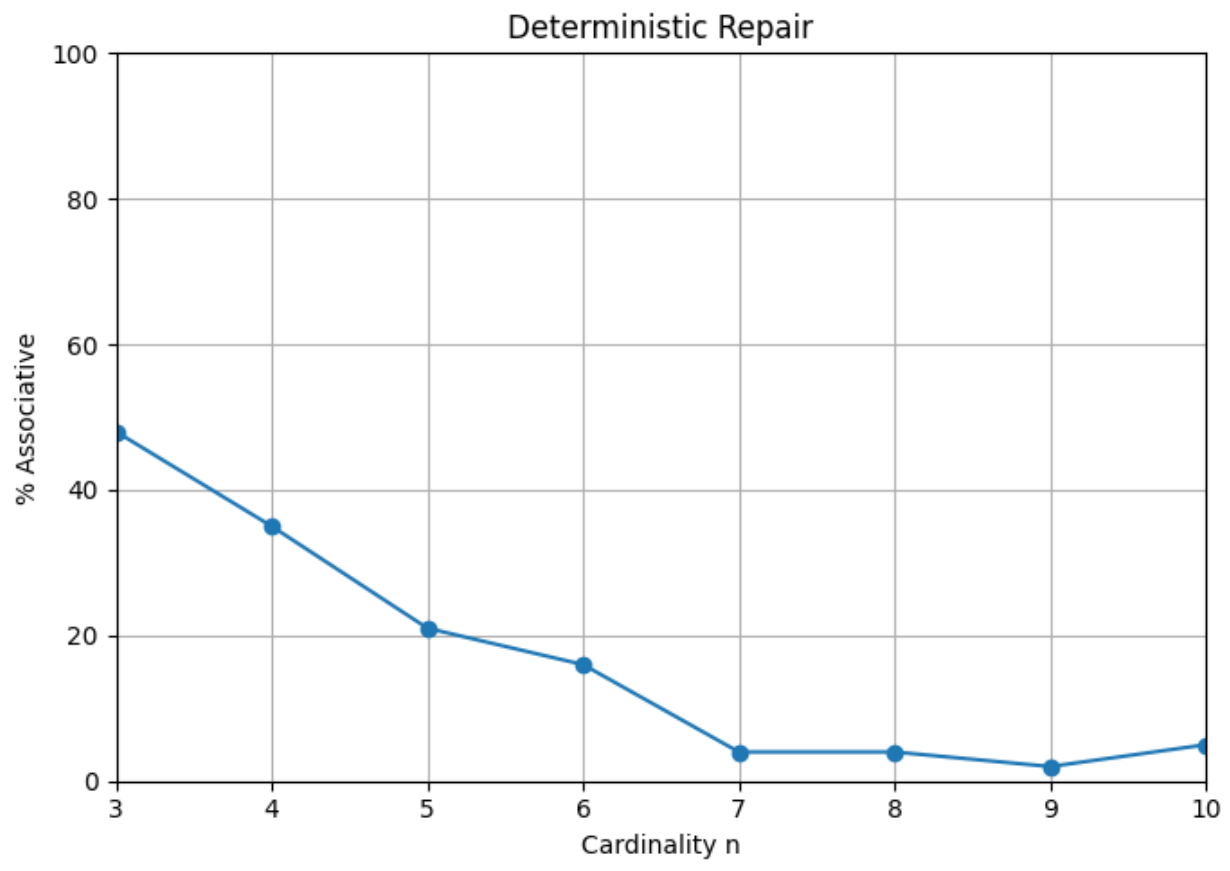}
        \caption{\centering Deterministic Healing (p=15\%)}
        \label{fig:assoc_vs_corr}
    \end{minipage}
    \hfill
    \begin{minipage}[b]{0.48\textwidth}
        \centering
        \includegraphics[width=\linewidth]{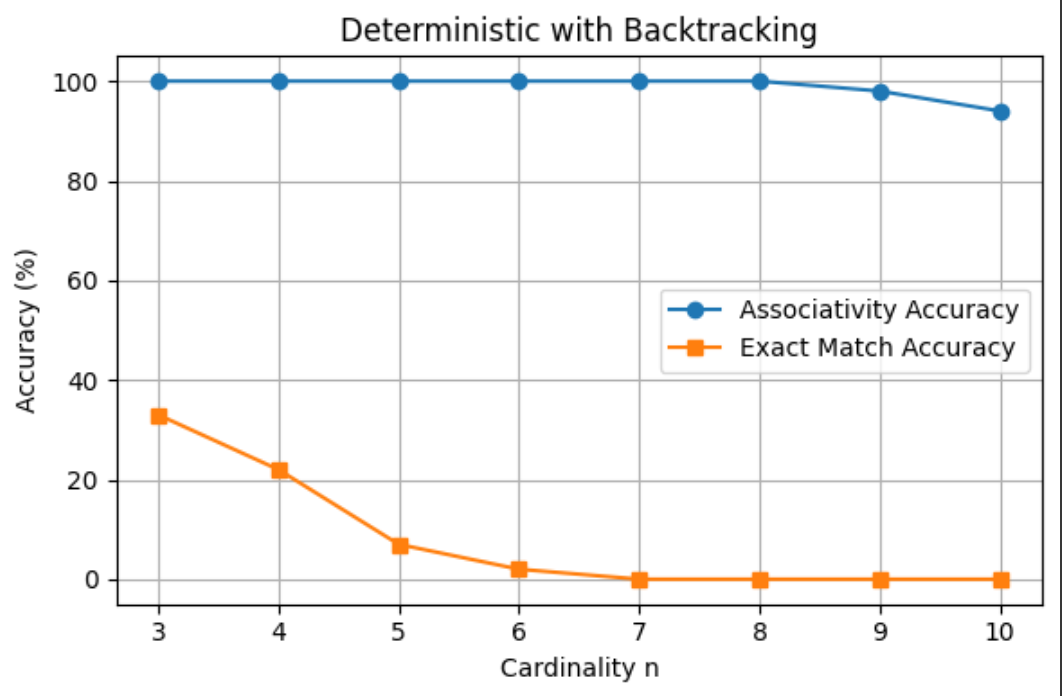}
        \caption{\centering Deterministic with Backtracking (p=15\%)}
        \label{fig:second_fig}
    \end{minipage}
\end{figure}

Figure \ref{fig:assoc_vs_corr} illustrates that while the deterministic repair strategy shows some efficiency at small cardinalities ($47\%$ at $n=3$), its performance deteriorates rapidly with increasing cardinality. By $n=7$, fewer than 5\% of the tables remain associative after repair and by $n=9$, the method is non-functional.

We also tested a backtracking variant (Figure \ref{fig:second_fig}): whenever a repair step failed, the algorithm retraced to the previous cell and tried the next candidate. This improved associativity rates in some cases, but at a cost: the repaired tables often diverged significantly from the ground truth and frequently collapsed into trivial semigroups (e.g., all entries equal). Thus, backtracking increased associativity but undermined fidelity, limiting its usefulness.

\section{Random Forest Classifier}

\begin{figure}
    \centering
    \includegraphics[width=1\linewidth]{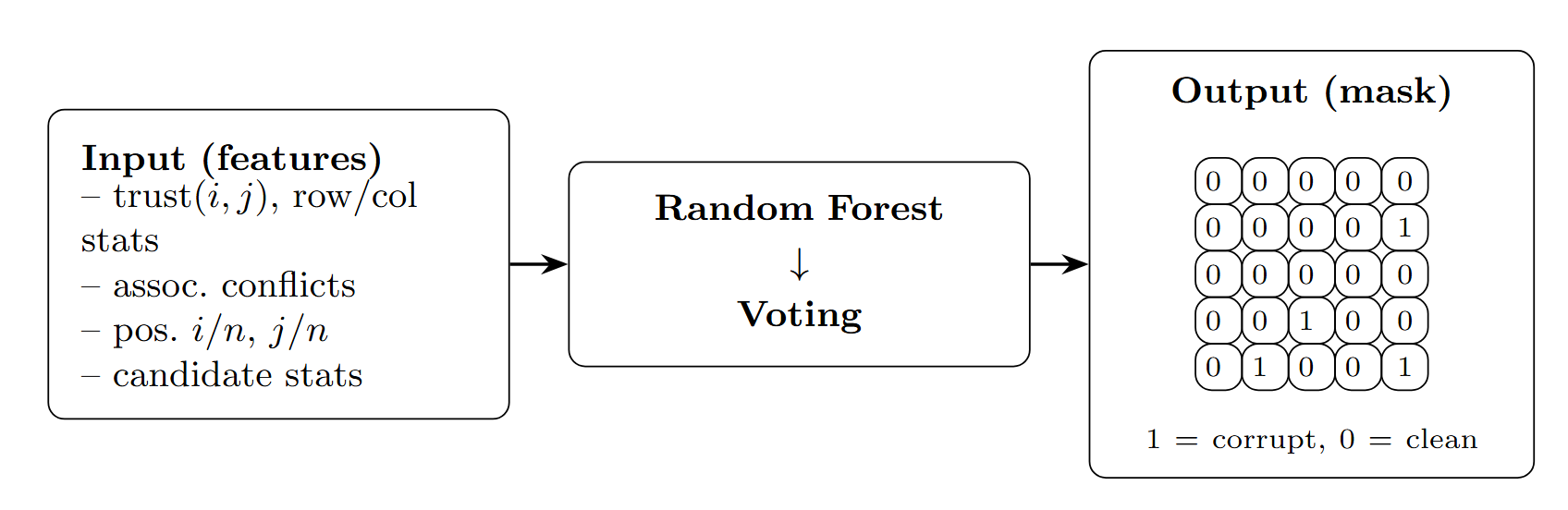}
    \caption{Random Forest Model Flow}
    \label{fig:rf-three-stage}
\end{figure}

\paragraph{Random Forest (RF) Algorithm.}
A Random Forest classifier is a ensemble of decision trees \cite{breiman2001randomforests,w3schoolsRF}, where each "tree" is trained on a random subset of the training data and a random subset of the features. 

In the language of machine learning, a \emph{feature} refers to a measurable attribute of each data point. Mathematically, if we represent an input as a vector $x \in \mathbb{R}^{d}$ in a $d$-dimensional feature space, then each coordinate $x_i$ is one feature. In this project, features include values such as trust scores, row/column indices, or local structural properties of the Cayley table.

Each tree learns a sequence of threshold tests (i.e., if trust $<0.3$ and row index $>5$, predict corrupted), and the ensemble combines their votes based off majority. Ultimately, it "classifies" each cell in the semigroup as corrupted or clean. The cells that are marked as corrupted are masked for repair.

\paragraph{RF prevents overfitting.} RF does this through two randomization techniques:
\begin{enumerate}
    \item Bootstrap Sampling: The data is randomly sampled so each tree sees a slightly different subset of the data. This reduces variance and causes less overfitting to noise in a single training set.
    \item Random Feature Selection: The RF uses a random subset of features when building each decision tree, adding another layer of randomness. This prevents any single feature from dominating. 
\end{enumerate}

\paragraph{Training.}
The RF uses this large collection of decision trees as its method of "training." It also uses tree splitting, finding the best way to split the data based on a criterion (i.e. Gini impurity or information gain (entropy)) Importantly, the shape/structure of each individual decision tree is not fixed and instead grown dynamically during training. In this way, the forest collectively encodes a diverse set of decision rules rather than relying on a single rigid model.  

Beyond the generic feature construction, some features are designed to capture higher-level structural information about the table. Row and column trust scores are crucial because if one cell is corrupted it causes multiple associative violations in its row and column.
For example, if T[2,3] is wrong, then for every k both (T[2,3], k) and (2, T[3,k]) are now mismatched. By including row/col trust as features, the RF can exploit these "global" signals. For instance, if an entire row $r$ is flagged as low-trust (many associativity checks fail when $r$ participates), then any cell $(r, j)$ or $(i, r)$ should be treated with higher suspicion.

\paragraph{Human-like Reasoning.}
Each tree in the ensemble can be viewed as a chain of "if-then" statements, mimicking a simple form of human reasoning. For instance, one tree may learn rules such as "if trust$<0.5 \rightarrow$ suspicious, if trust$<0.3 \rightarrow $ corrupted, else $\rightarrow$ clean." Another tree will discover different thresholds and feature combinations. While any individual tree may be noisy or overly specific, the aggregation of many such trees through majority voting produces robust predictions.  

\section{Hybrid Framework}

\paragraph{The Idea.} We now present our Hybrid framework, which integrates a Random Forest (RF) with a deterministic repair algorithm. The RF is used to predict which cells that are likely corrupted, while the deterministic procedure enforces associativity and performs the actual repair. In addition, the RF guides the deterministic step by indicating (i) which cells should be targeted for correction, (ii) which cells can be trusted and reused in the repair process, and (iii) how candidate repairs should be prioritized.
\begin{figure}[H]
\centering
\small
\resizebox{\linewidth}{!}{%
\begin{tikzpicture}[
  >=Stealth,
  node distance=6mm and 8mm,
  stage/.style={draw, rounded corners, align=center, fill=white,
                inner sep=6pt, font=\small, minimum width=2cm, minimum height=10mm},
  lab/.style={align=center, font=\footnotesize, text width=2cm}
]
\node[stage] (a) {\textbf{1.\ Generating}\\Data};
\node[stage, right=of a] (b) {\textbf{2.\ Training}\\Random Forest (RF)};
\node[stage, right=of b] (c) {\textbf{3.\ Deterministic}\\Repair};
\node[stage, right=of c] (d) {\textbf{4.\ Execution}};
\node[stage, right=of d] (e) {\textbf{5.\ Caching}};

\draw[->, thick] (a) -- (b);
\draw[->, thick] (b) -- (c);
\draw[->, thick] (c) -- (d);
\draw[->, thick] (d) -- (e);

\node[lab, below=1.5mm of a] {clean $\to$ corrupt, trust};
\node[lab, below=1.5mm of b] {RF threshold $\tau$};
\node[lab, below=1.5mm of c] {mask+fill, candidates};
\node[lab, below=1.5mm of d] {ML, subsemigroups, merge};
\node[lab, below=1.5mm of e] {save results, metrics};
\end{tikzpicture}
}%
\caption{Hybrid pipeline}
\label{fig:hybrid-flow}
\end{figure}
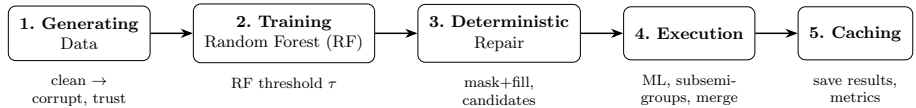

\paragraph{Generating Data.}
The first step in our framework is constructing the dataset of semigroups using Mace4 (see Section 3.1). We generate pairs of corrupted and clean semigroups. For every corrupted table, we also compute a \emph{trust map} (Section 3.2), which represents the reliability of each cell based on associativity violations. 
\paragraph{Training RF Classifier.}
The corrupted-clean pairs are split into test and train data. Feature vectors are extracted from each cell, consisting of the trust score, row/column indices, and candidate values. These features are passed into the Random Forest (RF) classifier. The RF predicts the probability of corruption and applies a threshold $\tau$ to mask low-trust cells by setting them to $-1$. This masking step focuses subsequent repair only on the most uncertain entries. 
\paragraph{Deterministic Repair.}
Masked tables are then passed through a deterministic repair procedure (Section 4). In this phase only the cells previously set to $-1$ are healed. The repair proceeds by iterating over candidate values consistent with associativity; if no valid candidate is found, the entry remains at $-1$. After the repair, we recompute the trust map for the partially healed table, creating an updated representation of reliability for use in the next stage.
\paragraph{Execution.}
At this stage, we incorporate the subsemigroup structure. Closure sets of the form  $s=\{i,j,k,i \cdot j,j \cdot k,(i \cdot j)\cdot k,i \cdot (j \cdot k)\}$ are generated (Section~3.3). A second healing pass is applied: within each subsemigroup, RF-derived weights are used to resolve conflicts, overwriting lower-weight values with higher-weight predictions. Finally, the full semigroup is reconstructed with evidence from all subsemigroups. For each product $i \cdot j$, the algorithm gathers candidate values across subgroups and assigns weights according to
\[
w = p(\text{correct}) \times \frac{1}{|s|} \times \text{trust},
\]
where $p(\text{correct})$ is given by the RF model, $|s|$ is the subgroup size, and the trust score reflects associativity consistency. The candidate with the highest weight is chosen as the final entry for the global table.  

\paragraph{Caching.}
Each experimental run is logged, recording the semigroup size $n$, corruption percentage $p\%$, and repair method used. These cached results allow for efficient post-hoc analysis: the tables do not need to be regenerated, and performance metrics can be easily produced and compared across methods. 

\section{Results}

\begin{figure}[H]
\centering
\begin{tikzpicture}
\begin{axis}[
    width=10.8cm,
    height=7.2cm,
    xlabel={Cardinality $n$},
    ylabel={Percentage (\%)},
    ymin=0, ymax=105,
    xtick={3,4,5,6,7,8,9,10,15,20},
    xmin=3, xmax=20,
    legend style={at={(0.5,-0.2)},anchor=north,legend columns=2},
    grid=both,
    grid style={dashed,gray!30}
]

\addplot+[orange, mark=circle*, ultra thick] coordinates {
(3,95) (4,95) (5,95) (6,95) (7,70) (8,80) (9,60) (10,60) (15,30) (20,30)
};
\addlegendentry{Hybrid: \% Fully Associative}

\addplot+[orange, mark=square*, dashed, thick] coordinates {
(3,98.89) (4,99.84) (5,99.92) (6,99.86) (7,98.63) (8,98.46) (9,99.14) (10,97.85) (15,82.51) (20,91.74)
};
\addlegendentry{Hybrid: Mean Associativity Fraction}

\addplot+[orange, mark=triangle*, dotted, thick] coordinates {
(3,93.89) (4,92.50) (5,83.00) (6,76.11) (7,81.53) (8,75.00) (9,72.22) (10,67.25) (15,53.13) (20,49.05)
};
\addlegendentry{Hybrid: Mean Cell Acc}

\addplot+[blue, mark=circle*, ultra thick] coordinates {
(3,47) (4,35) (5,20) (6,15) (7,5) (8,3) (9,2) (10,1) (15,0) (20,0)
};
\addlegendentry{Deterministic: \% Fully Associative}

\end{axis}
\end{tikzpicture}
\caption{Hybrid vs deterministic healing performance across cardinalities. Hybrid results include three metrics, while deterministic repair quickly collapses and is shown only by \% fully associative.}
\label{fig:hybrid-vs-det}
\end{figure}
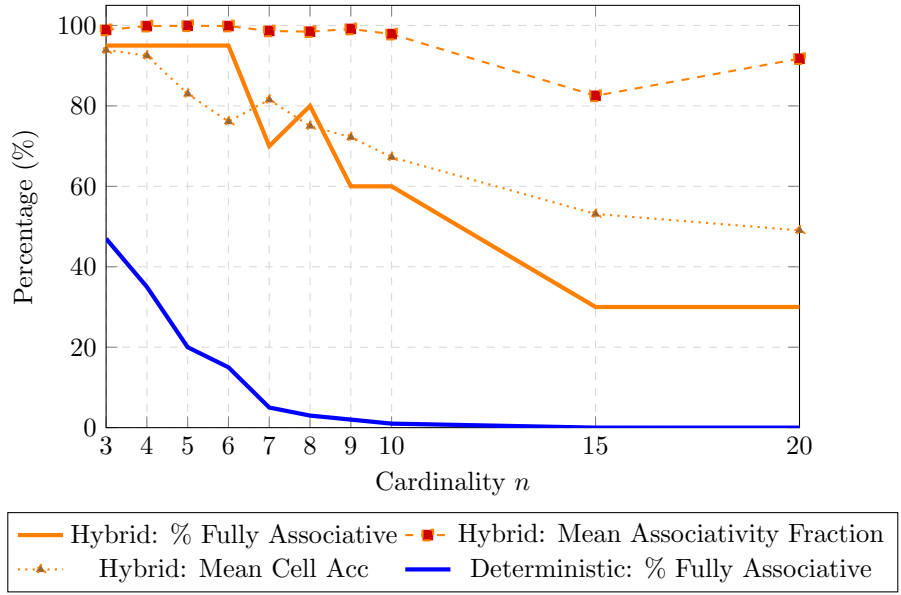

\paragraph{Evaluation.}
The hybrid pipeline substantially outperforms deterministic repair. For $n \leq 6$, more than $95\%$ of corrupted tables were fully repaired. Accuracy declines at higher $n$, but even at $n=10$ the hybrid approach achieved $60\%$ fully associative tables, compared to only $2\%$ for deterministic repair. This demonstrates that the combination of trust maps, Random Forest masking, and local subsemigroup healing is a better method of healing compared to plain deterministic.

\paragraph{Metrics for Success.}
Although the percentage of tables that are fully associative is the primary measure of success, we also consider two complementary metrics: (i) the mean per-cell accuracy, and the (ii) the mean associative fraction (the proportion of triples $(i,j,k)$ satisfying associativity). Despite the visible decline in fully associative rates for $n > 8$, both mean associativity fraction and mean per-cell accuracy remain remarkably high.

The mean per-cell accuracy is important to ensure the healed semigroups are fairly close to the original table, and prevents the algorithm from always producing the same or trivial semigroups. The associative fraction presents a more nuanced view on how effective the healing is even when tables aren't fully healed. For example, at $n$=20 the associative fraction is still high at 92\%. This shows that even when global associativity fails, the hybrid approach produces partially healed tables that preserve much of the local structure.

\begin{figure}[H]
\centering
\begin{tikzpicture}
\begin{axis}[
    width=0.85\linewidth,
    height=0.5\linewidth,
    xlabel={Healing Passes},
    ylabel={\% Associative Tables},
    symbolic x coords={Baseline, Healing~Pass~1, Healing~Pass~2},
    xtick=data,
    ymin=0, ymax=100,
    legend pos=south east,
    legend style={at={(0.98,0.5)}, anchor=west},
    grid=both,
    ymajorgrids=true,
    grid style=dashed,
]
\addplot coordinates {(Baseline,35) (Healing~Pass~1,95) (Healing~Pass~2,95)};
\addlegendentry{$n=3$}
\addplot coordinates {(Baseline,10) (Healing~Pass~1,90) (Healing~Pass~2,95)};
\addlegendentry{$n=4$}
\addplot coordinates {(Baseline,5) (Healing~Pass~1,85) (Healing~Pass~2,95)};
\addlegendentry{$n=5$}
\addplot coordinates {(Baseline,5) (Healing~Pass~1,85) (Healing~Pass~2,95)};
\addlegendentry{$n=6$}
\addplot coordinates {(Baseline,5) (Healing~Pass~1,60) (Healing~Pass~2,75)};
\addlegendentry{$n=7$}
\addplot coordinates {(Baseline,5) (Healing~Pass~1,65) (Healing~Pass~2,80)};
\addlegendentry{$n=8$}
\addplot coordinates {(Baseline,5) (Healing~Pass~1,20) (Healing~Pass~2,60)};
\addlegendentry{$n=10$}
\addplot coordinates {(Baseline,5) (Healing~Pass~1,15) (Healing~Pass~2,30)};
\addlegendentry{$n=15$}
\addplot coordinates {(Baseline,5) (Healing~Pass~1,15) (Healing~Pass~2,30)};
\addlegendentry{$n=20$}
\end{axis}
\end{tikzpicture}
\caption{Associativity recovery rates across healing passes for varying cardinalities.}
\label{fig:healing-comparison}
\end{figure}
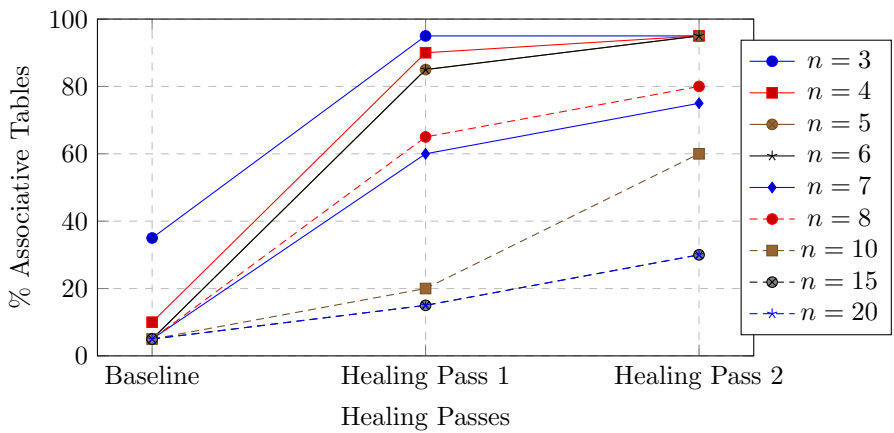

\paragraph{Effectiveness of each Healing Pass.}
We next evaluate the contribution of individual healing passes. Figure \ref{fig:healing-comparison} shows the percentage of fully associative tables after each pass. Healing Pass 1 (includes masking with RF and deterministic repair) achieves strong performance at small cardinality, but drops sharply at larger $n$ (e.g., 85\% at $n=5$ and only 20\% at $n=10$). By adding Healing Pass 2 (subsemigroups decomposition and reconstruction), performance improves substantially, reaching 95\% at $n=5$ and 60\% at $n=10$. This confirms that the multi-stage hybrid pipeline is necessary: neither pass alone is sufficient at larger scales, but their combination produces reliable results. 

\section{Discussion}
\subsection{Statistical Analysis}

In principle, under a naive uniformity assumption, one might expect the correct value for $i \cdot j$ to dominate with probability close to $1-\tfrac{1}{n}$. If this were true, deterministic repair should perform fairly well, improving as cardinality increases. However, as our results show, the method performs poorly even at moderate corruption levels and across all cardinalities.

To investigate, we asked whether one could predict the exact point where deterministic repair "dips" in accuracy. We framed this as the probability that an incorrect value receives at least as many votes as the true one. Formally, this is the "Exceeds $C$" probability, where $C=(1-p)n$ is the expected number of correct votes under corruption rate $p$. Using the binomial formula,
\[
\Pr[X \ge C] = \sum_{k=C}^{n} \binom{n}{k} 
\left(\tfrac{1}{n}\right)^{k}
\left(1-\tfrac{1}{n}\right)^{n-k},
\]
we find that the probability is astronomically small--for example, for $n=10$ and $p=0.15$, it is only $9.1 \times 10^{-9}$. In other words, purely random fluctuations cannot explain the collapse of deterministic repair.

The contradiction arises precisely because the \emph{uniformity assumption} is false. The calculation assumes that all $n$ values appear with equal frequency across the semigroup table. However, real semigroups are highly skewed—some elements occur disproportionately often while others appear rarely or not at all. Table \label{tab:freq-full} illustrates this, showing the large standard deviations of element frequencies compared to the uniform baseline.  

\begin{table}[H]
\centering
\begin{tabular}{c c c c}
\toprule
Value $v$ & Count $s_v$ & Frequency $f_v$ & Deviation $(f_v - 1/n)$ \\
\midrule
0 & 5  & 0.05 & -0.05 \\
1 & 32 & 0.32 & +0.22 \\
2 & 4  & 0.04 & -0.06 \\
3 & 4  & 0.04 & -0.06 \\
4 & 26 & 0.26 & +0.16 \\
5 & 27 & 0.27 & +0.17 \\
6 & 0  & 0.00 & -0.10 \\
7 & 0  & 0.00 & -0.10 \\
8 & 0  & 0.00 & -0.10 \\
9 & 2  & 0.02 & -0.08 \\
\bottomrule
\end{tabular}
\caption{Observed value frequencies for $n=10$ in one semigroup table. 
The uniform target is $1/n = 0.1$. Large deviations show that some values dominate 
while others never appear, breaking the uniformity assumption.}
\label{tab:freq-full}
\end{table}

\paragraph{Outcome.}
This analysis yields two conclusions:  
(1) deterministic repair cannot be predicted to fail at a single clean threshold, since real semigroups are far from uniform; and  
(2) the method is inherently unreliable, motivating the need for machine learning to handle non-uniform corruption patterns.  

\subsection{Value of Machine Learning}
We have analyzed the limitations of deterministic repair through statistical analysis but we must understand the ways in which machine learning provides a decisive advantage? We identify two key factors. 

\begin{enumerate}
    \item \textbf{Decision Making.} Deterministic Repair requires iterating exhaustively over every $(i,j,k)$ triple and applying associative checks. In contrast, the Random Forest Classifier learns a compact representation of the corruption patterns. This allows the algorithm to target only those cells predicted to be unreliable, avoiding the combinatorial cost of brute-force iteration.
    \item \textbf{Capturing Semigroup Complexity.} Semigroups exhibit highly non-uniform distributions: some values appear disproportionately often, while others may not appear at all. Deterministic repair assumes uniformity and applies fixed cutoffs, which fails in these cases. Machine Learning, by contrast, can combine trust scores, positional features, and candidate distributions to form an effective "summary" of the semigroup's structure. This representation enables robust predictions even when value frequencies are skewed. Similar multimodal approaches highlight the benefit of combining structural priors with learning \cite{koehler2024multimodal}.
\end{enumerate}

\paragraph{Design Choice.}

The sequence of deterministic repair followed by machine learning was chosen deliberately. Deterministic reasoning resolves trivial cases (e.g., zeros, idempotents, small closure sets) and reduces noise before more nuanced predictions are required. The Random Forest then addresses the harder cases, combining trust scores, positional features, and subgroup overlaps to generalize beyond fixed thresholds. Applied, the two methods complement each other and achieve performance that neither could reach individually.  

\paragraph{Future Directions.}In future experiments, we plan to extend the study to higher cardinality ($n>30$) and systematically vary corruption levels (beyond $p=15\%$). These extensions will test the scalability and robustness of the hybrid approach beyond the settings considered here.

\section*{Acknowledgements}
The author would like to thank Jasper van de Kreeke for supervision and guidance throughout this project. This work was supported by NWO Rubicon grant 019.232EN.029.

The author used a generative AI language tool for assistance in refining wording and preparing some figures. All aspects of the work remain the author’s own.

\newpage
\bibliographystyle{plain}
\bibliography{refs}

\end{document}